\documentclass[11pt]{article}

\usepackage{latexsym}
\usepackage{amsfonts}
\usepackage{enumerate}
\usepackage{multicol}
\usepackage{graphicx}
\usepackage{amssymb}
\usepackage{amsmath}
\usepackage{amsthm}
\usepackage{epic}
\usepackage{sidecap}
\usepackage{tikz}

\begin{document}
\parindent=0.2in
\parskip 0in

\begin{flushright}

{\huge {\bf Report on the 50th Annual USA Mathematical Olympiad}}

\vspace*{.2in}

{\Large B\'ELA BAJNOK} \\

{\small Gettysburg College \\  Gettysburg, PA 17325 \\ bbajnok@gettysburg.edu}

\end{flushright}

\vspace*{.2in}

The USA Mathematical Olympiad (USAMO) is the final round in the American Mathematics Competitions series for high school students, organized each year by the Mathematical Association of America.  The competition follows the style of the International Mathematics Olympiad: it consists of three problems each on two consecutive days, with an allowed time of four and a half hours both days.  

The 50th annual USAMO was given on Tuesday, April 13, 2021 and Wednesday, April 14, 2021.  This year, 288 students were invited to take the USAMO and, as in 2020, the competition was administered online.  The names of winners and honorable mentions, as well as more information on the American Mathematics Competitions program, can be found on the site https://www.maa.org/math-competitions.  Below we present the problems and solutions of the competition; a similar article for the USA Junior Mathematical Olympiad (USAJMO), offered to students in grade 10 or below, can be found in the January 2022 issue of the {\em College Mathematics Journal}.     

The problems of the USAMO are chosen -- from a large collection of proposals submitted for this purpose -- by the USAMO/USAJMO Editorial Board, whose co-editors-in-chief this year were Evan Chen and Jennifer Iglesias, with associate editors Ankan Bhattacharya, John Berman, Zuming Feng, Sherry Gong, Alison Miller, Maria Monks Gillespie, and Alex Zhai.  This year's problems were created by Ankan Bhattacharya, Mohsen Jamaali, Shaunak Kishore, Carl Schildkraut, Zoran Sunic, and Alex Zhai.  

The solutions presented here are those of the present author, relying in part on the submissions of the problem authors and members of the editorial board.
Each problem was worth 7 points; the nine-tuple $(n; a_7, a_6, a_5, a_4, a_3, a_2, a_1, a_0)$ states the number of students who submitted a paper for the relevant problem, followed by the numbers who scored $7, 6, \dots, 0$ points, respectively.

\vspace{.1in}

\noindent {\bf Problem 1} $(226; 140, 0, 0, 0, 1, 0, 11, 74)$; {\em proposed by Ankan Bhattacharya}.     Rectangles $BCC_1B_2$, $CAA_1C_2$, and $ABB_1A_2$ are erected outside an acute triangle $ABC$. Suppose that
	\[ \angle BC_1C + \angle CA_1A + \angle AB_1B = 180^\circ. \]
Prove that lines $B_1C_2$, $C_1A_2$, and $A_1B_2$ are concurrent.

\vspace{.1in}

\noindent {\em Solution}.  Let $\omega_A$, $\omega_B$, and $\omega_C$ be the circumcircles of rectangles $BCC_1B_2$, $CAA_1C_2$, and $ABB_1A_2$, respectively.  Define $P$ to be the foot of the altitude from $A$ to $B_1C_2$.

\begin{figure}[ht]
\centering
\includegraphics{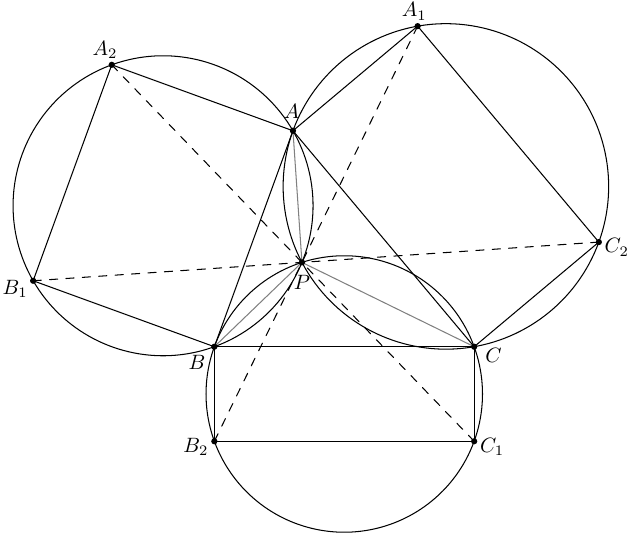}
\caption{The main circles used in Problem 1}
\end{figure}

Observe that, since $\angle ACC_2$ and $\angle APC_2$ are both $90^\circ$, $P$ lies on $\omega_B$ by the inscribed angle theorem and, therefore, $\angle APC$ and $\angle CA_1A $ are supplementary angles.  A similar argument shows that $P$ lies on $\omega_C$ and thus $ \angle APB$ and $\angle AB_1B$ are supplementary angles as well.

But then 
\[ \angle BPC = 360^\circ - (\angle APC + \angle APB) =  \angle CA_1A + \angle AB_1B,\]
which, by the given equation, yields 
$$ \angle BPC = 180^\circ - \angle BC_1C, $$ and thus $P$ lies on $\omega_A$ as well.  Therefore, $P$ is the (unique) common point of all three circles.

Similar arguments would prove that the feet of the altitudes from $B$ and $C$ to $C_1A_2$ and $A_1B_2$, respectively, are on each of the three circles, and thus must coincide with $P$.  But then lines $B_1C_2$, $C_1A_2$, and $A_1B_2$ are concurrent, as claimed.   

\vspace{.1in}

\noindent {\bf Problem 2} $(199; 76, 9, 1, 4, 13, 18, 1, 77)$; {\em proposed by Zoran Sunic}.   The Planar National Park is a subset of the Euclidean plane
	consisting of several trails which meet at junctions.
	Every trail has its two endpoints at two different junctions,
	whereas each junction is the endpoint of exactly three trails.
	Trails only intersect at junctions
	(in particular, trails only meet at endpoints).
	Finally, no trails begin and end at the same two junctions.

\begin{figure}[ht]
	\centering
	\includegraphics[scale=2]{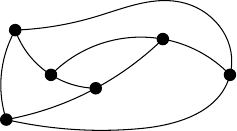}
\caption{An example of one possible layout of the park is shown
	in which there are six junctions and nine trails.}
\label{fig M2}
\end{figure}

	A visitor walks through the park as follows:
	she begins at a junction and starts walking along a trail.
	At the end of that first trail, she enters a junction and turns left.
	On the next junction she turns right, and so on,
	alternating left and right turns at each junction.
	She does this until she gets back to the junction where she started.
	What is the largest possible number of times she could have entered
	any junction during her walk, over all possible layouts of the park?

\vspace{.1in}

\begin{figure}[ht]
	\centering
	\includegraphics[width=0.4\textwidth]{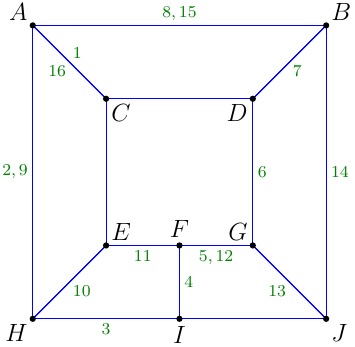}
	\caption{An example achieving three visits}
	\label{fig M2 construction}
\end{figure}

\noindent {\em Solution}.  The answer is three times.
We begin by exhibiting an example of a park layout which features three visits.
Sketched in Figure~\ref{fig M2 construction} is one of many possible constructions.
The path starts from $C$ and walks toward $A$,
and continues as follows:
\[ C \to A \to H \to I \to F \to G \to D \to B \to A \to H \to E \to F \to G
	\to J \to B \to A \to C \]
As we see, this path visits $A$ three times.

We will prove that the visitor cannot visit any junction more than three times.  (This trivially holds fore the initial/terminal junction.)  Note that if a junction were to be visited four times or more, then this would mean four or more arrivals and four or more departures, which is only possible if at least one of the three trails meeting at that junction would have had to be traversed at least three times.  Therefore, it suffices to show that no trail can be on the visitor's path more than twice.

Suppose, indirectly, that there is a trail that the visitor walked on three or more times.  (Note that this trail cannot be adjacent to the junction where her walk started.)  This then implies that at least two of those times she turned in the same direction (left or right) when she reached the end of the trail.  Let's assume then that her $m$-th and $n$-th trail during her walk is the same for some $2 \leq m<n$ with the same turn at the end; we may further assume that this is the trail with the smallest possible $m$.  Let $A$ and $B$ denote the junctions at the two ends of this trail.       

Now if she walked along this trail both times in the same direction, say from $A$ to $B$, and made the same turn at the end (e.g., left), then her $(m-1)$-st and $(n-1)$-st trails were also the same, and she made the same turn when she got to $A$ (right).  This contradicts the minimality of $m$.  On the other hand, if once she walked from $A$ to $B$ and then later from $B$ to $A$, but turning in the same direction at the end both times (e.g., left), then her $(m-1)$-st trail was the same as her $(n+1)$-st, and she made the same turn at the ends of these two trails as well (right), again 
contradicting the minimality of $m$.  This completes our proof.

\vspace{.1in}

\noindent {\bf Problem 3} $(179; 7, 4, 1, 0, 1, 4, 107, 55)$;  {\em proposed by Shaunak Kishore and Alex Zhai}.   Let $n \ge 2$ be an integer.
	An $n \times n$ board is initially empty.
	Each minute, you may perform one of three moves:
\begin{itemize}
	\item If there is an L-shaped tromino region
		of three cells without stones on the board
		(see Figure \ref{fig M3 problem}; rotations not allowed),
		you may place a stone in each of those cells.
\begin{figure}[ht]
	\centering
	\includegraphics{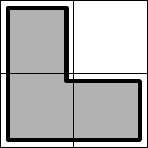}
	\caption{Three cells forming an L-shaped tromino}
	\label{fig M3 problem}
\end{figure}
	\item If all cells in a column have a stone, you may remove all stones from that column.
	\item If all cells in a row have a stone, you may remove all stones from that row.
	\end{itemize}

For which $n$ is it possible that, after some nonzero number of
	moves, the board has no stones?

\vspace{.1in}

\noindent {\em Solution.}  We claim that the answer is all positive integers $n$ that are divisible by 3.  First, we show that the procedure is possible in each of these cases.

\begin{figure}[ht]
	\centering
	\includegraphics{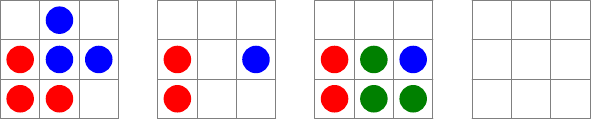}
	\caption{The four-step procedure that clears all stones}
 \label{fig-mo3}
\end{figure}

When $n$ is divisible by $3$,
one may divide the board into $3 \times 3$ sub-squares; 
for brevity, let us refer to these $(n/3)^2$ sub-squares as \emph{cages}.
We then follow the procedure illustrated in Figure \ref{fig-mo3}, as follows.
\begin{itemize}
	\item First, we put two non-overlapping L-trominoes in each cage, as shown in the first step.
	\item This causes every center column of each cell to be completely filled.
		Thus, we may remove all $n/3$ columns which correspond
		to the center columns of cages, as shown in the second step.
	\item In each cage, we then place one L-tromino as shown in the third step.
	\item Now the board consists of $2n/3$ completely filled rows,
		so we may eliminate them all.
\end{itemize}

We now prove that if after some sequence of moves no stones remain, then $n$ must be a multiple of 3.  We will employ what is usually called the {\em polynomial method}; in particular, we make use of the following.

{\bf Lemma}.  {\em Consider the polynomial $$f(x,y)=\sum_{i=0}^{n_1} \sum_{j=0}^{n_2} d_{i,j} x^iy^j,$$ where the coefficients $d_{i,j}$ are real numbers and $d_{n_1,n_2} \neq 0$.  If $A_1$ and $A_2$ are sets of real numbers with $|A_1| > n_1$ and $|A_2| > n_2$, then there are elements $a_1 \in A_1$ and $a_2 \in A_2$ for which $f(a_1,a_2) \neq 0$.}\footnote{Note that this lemma is the two-variable version of the well-known fact that a nonzero polynomial cannot have more roots than its degree.  For a simple proof and a variety of applications see, for example, Section 12.3 in  \cite{FreGya:2020a}.}

\vspace*{0.1in} 

Let us introduce some notations.  We parametrize the cells of the board by letting $(i,j)$ denote the position of the cell in the $i$-th column (counting from the left) and the $j$-th row (counting from the bottom).  We then associate each state of the board with the polynomial $$A(x,y)=\sum_{i=1}^n \sum_{j=1}^n c_{i,j} x^{i-1}y^{j-1},$$where $c_{i,j}$ is 1 when there is a stone in cell $(i,j)$ and 0 otherwise.  This allows us to think of the chain of moves as the sequence $\left( A_k (x,y) \right)_{k=0}^m$ where $A_0 (x,y) =0$ (representing the initial position of the board), $A_m (x,y) =0$ (expressing the fact that there are no stones on the board after $m$ moves for some $m \in \mathbb{N}$), and where $A_k (x,y) $ results from $A_{k-1} (x,y) $ in one of the following ways:
\begin{itemize}
	\item $A_k (x,y) =A_{k-1} (x,y) +x^{i-1}y^{j-1}(1+x+y)$ if a tromino was placed on the board with its lower left corner at position $(i,j)$;
	\item $A_k (x,y) =A_{k-1} (x,y) -x^{i-1}(1+y+y^2+\cdots + y^{n-1})$ if all stones got removed from column $i$; and
	\item $A_k (x,y) =A_{k-1} (x,y) -y^{j-1}(1+x+x^2+\cdots + x^{n-1})$ if all stones got removed from row $j$.
\end{itemize}   
We need a few more notations. For $i=1,2,\dots, n-1$ and $j=1,2,\dots, n-1$, we let $a_{i,j}$ denote the number of times a tromino was added with its lower-left corner at position $(i,j)$; we then set $$P(x,y)=\sum_{i=1}^{n-1} \sum_{j=1}^{n-1} a_{i,j} x^{i-1}y^{j-1}.$$  Furthermore, we set $r(j)$ and $c(i)$ equal to the number of times the $j$-th row and $i$-th column were cleared, respectively, and define $Q(x)=\sum_{i=1}^n c_i x^{i-1}$ and $R(y)=\sum_{j=1}^n r_j y^{j-1}$.  With these notations, the fact that our procedure succeeded can be stated by the equation
$$P(x,y) (1+x+y) - Q(x)(1+y+y^2+\cdots + y^{n-1}) - R(y)(1+x+x^2+\cdots + x^{n-1})=0.$$

Take $A$ to be the set of $n$-th roots of unity other than $1$; that is, the $n-1$ distinct complex numbers $a$ for which $$\frac{a^n-1}{a-1}=1+a+a^2+\cdots+a^{n-1}=0.$$  Since $P(x,y)$ is a nonzero polynomial with $x$-degree and $y$-degree at most $n-2$, our lemma above guarantees elements $a_1, a_2 \in A$ for which    
$P(a_1,a_2) \neq 0$.  But since substituting $x=a_1$ and $y=a_2$ into our equation yields $$P(a_1,a_2)(1+a_1+a_2)=0,$$ this can only occur when $1+a_1+a_2=0$.  Therefore, the imaginary parts of $a_1$ and $a_2$ are negatives of one another; recalling that both numbers have norm 1, this then implies that their real parts have the same absolute value.  But these real parts must then both be negative, and in fact equal to $-\frac{1}{2}$, so $a_1$ and $a_2$ are $-\frac{1}{2} \pm \frac{\sqrt{3}}{2} \mathrm{i}$.  We thus got that $a_1$ and $a_2$ are third roots of unity.  However, this can only happen if $n$ is divisible by 3: indeed, if $n=3q+r$ for some integers $q$ and $r$ with $r \in \{0,1,2\}$, then $a_1^r=a_1^n/a_1^{3q}=1$, 
which is only possible when $r=0$, as claimed.

\vspace{.1in}

\noindent {\bf Problem 4} $(240; 121, 14, 5, 0, 1, 37, 13, 49)$; {\em proposed by Carl Schildkraut}.   A finite set $S$ of positive integers has the property that,
	for each $s\in S$, and each positive integer divisor $d$ of $s$,
	there exists a unique element $t\in S$ satisfying $\gcd(s,t) = d$.
	(The elements $s$ and $t$ could be equal.)

	Given this information, find all possible values for the
	number of elements of $S$.

\vspace{.1in} 

\noindent {\em Solution.}  We claim that the possible sizes are 0 and the nonnegative integer powers of 2.  Since $S=\emptyset$ and $S=\{1\}$ obviously work, we need to show that a set $S$ of size $n \geq 2$ satisfying the requirements exists if, and only if, $n=2^k$ for some positive integer $k$.   

We start by verifying that these values are indeed possible.  For a given positive integer $k$, we construct a set of size $2^k$ as follows.  Suppose that $p_1, q_1, p_2, q_2, \dots, p_k, q_k$ are $2k$ pairwise distinct positive primes; for an ordered pair of subsets $(I,J)$ of $[k]=\{1,2, \dots, k\}$, we will use the notation 
$$s(I,J)=\prod_{i \in I} p_i \cdot \prod_{j \in J} q_j.$$  (Recall that the empty product equals 1.)  
We then consider 
the set $$S=\{s(I,J) \mid I \subseteq [k], J=[k] \setminus I\}.$$  Since the elements of $S$ are then in a bijection with the subsets $I$ of $[k]$, we see that $|S|=2^k$.  We need to show that $S$ satisfies the required property.  

Given an element $s(I,J)$ of $S$, we see that its positive divisors are of the form $s(I_0,J_0)$ where $I_0 \subseteq I$ and $J_0 \subseteq J$.  Note also that, since $I$ and $J$ form a partition of $[k]$, the sets $I'=I_0 \cup (J \setminus J_0)$ and $J'=J_0 \cup (I \setminus I_0)$ do as well, and thus 
$s(I',J')$ is an element of $S$; in fact it is the unique element of $S$ whose greatest common divisor with $s(I,J)$ equals $s(I_0,J_0)$.   Therefore, the set $S$ we constructed satisfies the requirement of the problem.

It remains to be shown that if $S$ is a set satisfying the property and it has size $n  \geq 2$, then $n=2^k$ for some positive integer $k$.  Let $s \geq 2$ be any element of $S$, and let $p$ be any positive prime divisor of $s$.  We can then write $s=p^e \cdot u$ for some positive integers $e$ and $u$ where $u$ is not divisible by $p$.  We claim that $e=1$.  

Denoting by $d(m)$ the number of positive divisors of a positive integer $m$, we have $d(s)=(e+1) \cdot d(u)$; in fact, $s$ has exactly $e  \cdot d(u)$ positive divisors that are divisible by $p$ and $d(u)$ that are not.  By our assumption, there is an element $t$ of $S$ for which $\gcd (t,s)=p$.  Let us assume that $e \geq 2$; we can then see that $t=p \cdot v$ for some positive integer $v$ that is not divisible by $p$.  Furthermore, $d(t)=2d(v)$, and $t$ has exactly $d(v)$ positive divisors that are divisible by $p$ and also $d(v)$ that are not.

Now according to our requirement for $S$, the positive divisors of any element $s$ of $S$ are in a one-to-one correspondence with the elements of $S$ via the map $a \mapsto \gcd(a,s)$, and thus $|S|=d(s)$.  (In particular, all elements of $S$ must have the same number of positive divisors.)  Furthermore, for any prime divisor $p$ of $s$, we have $p|a$ if and only if $p| \gcd(a,s)$.  Therefore, $S$ has exactly $e  \cdot d(u)$ elements that are divisible by $p$ and $d(u)$ that are not.  With the same reasoning, $S$ has exactly $d(v)$ elements that are divisible by $p$ and $d(v)$ that are not.  But then $d(u)=d(v)$ and $e  \cdot  d(u)=d(v)$, from which $e=1$.  

This establishes the fact that each element of $S$ is a product of the same number of pairwise distinct prime numbers.  If $s \in S$ is the product of $k$ distinct primes, then $|S|=d(s)=2^k$.  This completes our proof.

\vspace{.1in}

\noindent {\bf Problem 5} $(215; 95, 10, 0, 1, 2, 29, 2, 76)$; {\em proposed by Mohsen Jamaali}.   Let $n \ge 4$ be an integer.
	Find all positive real solutions to the following
	system of $2n$ equations:
	\begin{align*}
		a_1 = &\frac{1}{a_{2n}} + \frac{1}{a_{2}}, & a_2 &= a_1 + a_3, \\[1ex]
		a_3 = &\frac{1}{a_{2}} + \frac{1}{a_{4}}, & a_4 &= a_3 + a_5, \\[1ex]
		a_5 = &\frac{1}{a_{4}} + \frac{1}{a_{6}}, & a_6 &= a_5 + a_7, \\[1ex]
& \vdots & & \vdots \\
		a_{2n-1} = &\frac{1}{a_{2n-2}} + \frac{1}{a_{2n}}, & a_{2n} &= a_{2n-1} + a_1.
	\end{align*}

\vspace{.1in}

\noindent {\em First solution.}  It is easy to see that a solution is provided by $a_1=a_3=\cdots=a_{2n-1}=1$ and $a_2=a_4=\cdots=a_{2n}=2$; we prove that there are no others.

Taking indices modulo $2n$ and eliminating terms of odd indices, for each $i=1,2,\dots,n$ we have
\begin{eqnarray}
a_{2i}=\frac{1}{a_{2i-2}}+\frac{2}{a_{2i}}+\frac{1}{a_{2i+2}}. 
\end{eqnarray}
Adding up these equations yields 
\begin{eqnarray}
\sum_{i=1}^n a_{2i}=\sum_{i=1}^n \frac{4}{a_{2i}}.
\end{eqnarray}
According to the harmonic mean--arithmetic mean inequality, 
\begin{eqnarray}
\frac{n}{\sum_{i=1}^n \frac{1}{a_{2i}}} \leq \frac{\sum_{i=1}^n a_{2i}}{n}, 
\end{eqnarray}
so by (2) we get $\sum_{i=1}^n a_{2i} \geq 2n$.

Now dividing both sides of (1) by $a_{2i}$ yields
\begin{eqnarray}
1=\frac{1}{a_{2i-2}a_{2i}}+\frac{2}{a_{2i}^2}+\frac{1}{a_{2i}a_{2i+2}},
\end{eqnarray}
and adding the equations results in
\begin{eqnarray}
n=\sum_{i=1}^n \left( \frac{1}{a_{2i}}+ \frac{1}{a_{2i+2}} \right)^2.
\end{eqnarray}
Now we use the quadratic mean--arithmetic mean inequality, which gives
\begin{eqnarray}
\frac{1}{n} \cdot \sum_{i=1}^n \left( \frac{1}{a_{2i}}+ \frac{1}{a_{2i+2}} \right)^2 \geq \left( \frac{\sum_{i=1}^n \left( \frac{1}{a_{2i}}+ \frac{1}{a_{2i+2}} \right)}{n} \right)^2, 
\end{eqnarray}
so by (5) and (2) we get
$$1 \geq \left( \frac{ \sum_{i=1}^n \frac{2}{a_{2i}}}{n} \right)^2 =  \left( \frac{ \frac{1}{2} \cdot  \sum_{i=1}^n a_{2i} }{n} \right)^2 $$
and thus $\sum_{i=1}^n a_{2i} \leq 2n$.

This means that each of our inequalities is an equality, and therefore $a_{2i}=2$ for all $i$.  This in turn implies that $a_{2i-1}=1$ for all $i$, as claimed.  

\vspace{.1in}

\noindent {\em Second solution.}  We write $m= \min \{a_{2i} \mid 1 \leq i \leq n\}$ and $M= \max \{a_{2i} \mid 1 \leq i \leq n\}$, and assume that $m=a_{2j}$ and $M=a_{2k}$.  Then equation (1) from the first solution yields 
$$m=\frac{1}{a_{2j-2}}+\frac{2}{m}+\frac{1}{a_{2j+2}} \geq \frac{1}{M}+\frac{2}{m}+\frac{1}{M}$$ and 
$$M=\frac{1}{a_{2k-2}}+\frac{2}{M}+\frac{1}{a_{2k+2}} \leq \frac{1}{m}+\frac{2}{M}+\frac{1}{m}.$$
Therefore, $$m \geq \frac{2}{m}+\frac{2}{M} \geq M,$$ which can only occur when $m=M$.  Therefore, all $a_{2i}$ are equal, from which  $a_2=a_4=\cdots=a_{2n}=2$ and $a_1=a_3=\cdots=a_{2n-1}=1$, as claimed.

\vspace{.1in}

\noindent {\bf Problem 6} $(133; 17, 2, 1, 1, 3, 1, 5, 103)$; {\em proposed by Ankan Bhattacharya}.    Let $ABCDEF$ be a convex hexagon satisfying
	$\overline{AB} \parallel \overline{DE}$,
	$\overline{BC} \parallel \overline{EF}$,
	$\overline{CD} \parallel \overline{FA}$, and
	\[ AB \cdot DE = BC \cdot EF = CD \cdot FA. \]

	Let $X$, $Y$, and $Z$ be the midpoints
	of $\overline{AD}$, $\overline{BE}$, and $\overline{CF}$.
	Prove that
	the circumcenter of $\triangle ACE$,
	the circumcenter of $\triangle BDF$, and
	the orthocenter of $\triangle XYZ$ are collinear.
\vspace{.1in}

\noindent {\em Solution.}  We will prove that the orthocenter of $\triangle XYZ$ is in fact the midpoint of the segment connecting the circumcenter of $\triangle ACE$ and 
	the circumcenter of $\triangle BDF$.

\begin{figure}[ht]
	\centering
	\includegraphics{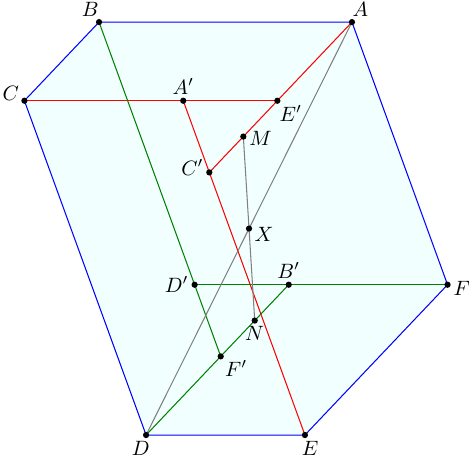}
	\caption{Parallelograms used in the solution to Problem 6}
	\label{fig M6}
\end{figure}

For each pair of adjacent sides of the hexagon, we construct a parallelogram with these two sides; this results in parallelograms $ABCE'$, $BCDF'$, $CDEA'$, $DEFB'$, $EFAC'$, and $FABD'$, as shown in Figure~\ref{fig M6}.  
(To aid legibility, Figure~\ref{fig M6} is intentionally \emph{not} drawn to scale.)  

Note that the assumption that $\triangle XYZ$ is nondegenerate implies that no two opposite sides of the hexagon have the same length; therefore, $\triangle A'C'E'$ and $\triangle B'D'F'$ are also nondegenerate triangles.   It is easy to see that these two triangles are translations of one another, as corresponding sides are parallel and have the same length: for example, $\overline{A'E'} \parallel \overline{D'B'}$ and, assuming without loss of generality that $AB > DE$, we have
$$A'E'=CE'-CA'=AB-DE=D'F-B'F=D'B'.$$  

Now let $M$ and $N$ be the midpoints of  $\overline{C'E'}$ and $\overline{B'F'}$, respectively.  Then $MAND$ is a parallelogram, because $\overline{AM}$ and $\overline{DN}$ are parallel, and  
$$AM=AE'+\frac{1}{2}E'C'=DF'+\frac{1}{2}F'B'=DN.$$  Therefore, $X$ (the midpoint of $\overline{AD}$) is the midpoint of $\overline{MN}$; similarly, $Y$ is the midpoint of the segment connecting the midpoint of $\overline{A'C'}$ and the midpoint of $\overline{F'D'}$, and $Z$ is the midpoint of the segment connecting the midpoint of $\overline{A'E'}$ and the midpoint of $\overline{B'D'}$.  Therefore, the orthocenter of $\triangle XYZ$ is the midpoint of the segment connecting the orthocenters of the medial triangles of $\triangle A'C'E'$ and $\triangle B'D'F'$.  Recall that the orthocenter of the medial triangle of a triangle $\triangle$ is the circumcenter of $\triangle$, hence the orthocenter of $\triangle XYZ$ is the midpoint of the segment connecting the circumcenters of $\triangle A'C'E'$ and $\triangle B'D'F'$.  We can thus complete our proof by showing that the circumcenters of $\triangle ACE$ and $\triangle A'C'E'$ coincide and that the circumcenters of $BDF$ and $B'D'F'$ coincide.  We show the first of these as the second claim can be done similarly. 

With $\omega$ denoting the circumcircle of  $\triangle A'C'E'$, set $r$ equal to the radius of $\omega$, and let $d_1$, $d_2$, and $d_3$ be the distances of $A$, $C$, and $E$ from the center of $\omega$, respectively. The power of $A$ to $\omega$ is then 
$$d_1^2-r^2=AE' \cdot AC' = BC \cdot EF;$$ similarly, we have 
$$d_2^2-r^2=CE' \cdot CA' = AB \cdot DE$$ and $$d_3^2-r^2=EA' \cdot EC' = CD \cdot FA.$$  According to our assumption, the three quantities are the same, which implies that $A$, $C$, and $E$ have the same distance from the center of $\omega$, and thus the circumcenters of $\triangle ACE$ and $\triangle A'C'E'$ coincide, as claimed.  This completes our proof.

\vspace*{.2in}

\noindent {\scriptsize {\bf Acknowledgments.} The author wishes to express his immense gratitude to everyone who contributed to the success of the competition: the students and their teachers, coaches, and proctors; the problem authors; the USAMO Editorial Board; the graders; the Art of Problem Solving; and the AMC Headquarters of the MAA.  I am also grateful to Evan Chen for proofreading this article and for providing the figures.}

\vspace*{.2in}

\vspace*{.2in}

\noindent {\bf B\'ELA BAJNOK} (MR Author ID: 314851) is a professor of mathematics at Gettysburg College and the director of the American Mathematics Competitions program of the MAA.

\end{document}